\documentclass[11pt]{article}

\usepackage{color}
\usepackage{latexsym}
\usepackage{dsfont}
\usepackage{amssymb}
\usepackage{graphicx}
\usepackage{amsmath, amsfonts,amssymb,theorem,euscript,array,enumerate,amsfonts,mathrsfs}
\usepackage{hyperref}

\newtheorem{Theorem}{Theorem}[part]

\newtheorem{Lemma}{Lemma}[part]

\newtheorem{Remark}{Remark}[part]

\def \trans{^{\scriptscriptstyle{\intercal}}}

\def \Int{\displaystyle\int}

\def \Inf{\displaystyle\inf}

\def \b1{\bf{1}}
\def \bx{\bf{x}}
\def \by{\bf{y}}

\def \trans{^{\scriptscriptstyle{\intercal }}}

\def \N{\mathbb{N}}
\def \R{\mathbb{R}}

\def \E{\mathbb{E}}
\def \F{\mathbb{F}}

\def \P{\mathbb{P}}

\def\esssup_#1{\underset{#1}{\mathrm{ess\,sup\, }}}

\def\argmin_#1{\underset{#1}{\mathrm{argmin\, }}}

\def \Ac{{\cal A}}
\def \Bc{{\cal B}}

\def \Fc{{\cal F}}

\def \Pc{{\cal P}}

\def \Nc{{\cal N}}

\def \eps{\varepsilon}

\def \ep{\hbox{ }\hfill$\Box$}

\def\reff#1{{\rm(\ref{#1})}}

\def\bx{{\bf x}}
\def\by{{\bf y}}

\def\beqs{\begin{eqnarray*}}
\def\enqs{\end{eqnarray*}}
\def\beq{\begin{eqnarray}}
\def\enq{\end{eqnarray}}

\begin{document}

 \title{Discrete time McKean-Vlasov control problem: \\ a dynamic programming approach}

\author{Huy\^en PHAM
\\\small  Laboratoire de Probabilit\'es et
 \\\small  Mod\`eles Al\'eatoires, CNRS, UMR 7599
 \\\small  Universit\'e Paris Diderot
 \\\small  pham at math.univ-paris-diderot.fr
\\\small  and CREST-ENSAE
\and
Xiaoli WEI 
\\\small  Laboratoire de Probabilit\'es et
 \\\small  Mod\`eles Al\'eatoires, CNRS, UMR 7599
 \\\small  Universit\'e Paris Diderot
 \\\small  tyswxl at gmail.com
}


\maketitle

\date{}

\begin{abstract}
We consider the stochastic optimal control problem of  nonlinear mean-field  systems in  discrete time.  We reformulate the problem into a deterministic control problem with marginal distribution as controlled state variable, 
and prove that dynamic programming principle holds in its general form. We apply our method for solving explicitly  the mean-variance portfolio selection and the multivariate linear-quadratic McKean-Vlasov control problem.
\end{abstract}

\vspace{5mm}


\vspace{5mm}



\vspace{5mm}

\section{Introduction}

\setcounter{equation}{0} \setcounter{Assumption}{0}
\setcounter{Theorem}{0} \setcounter{Proposition}{0}
\setcounter{Corollary}{0} \setcounter{Lemma}{0}
\setcounter{Definition}{0} \setcounter{Remark}{0}

The problem studied in this paper concerns the optimal control of nonlinear stochastic dynamical systems in discrete time of McKean-Vlasov type.  
Such topic is related to the modeling of  collective behaviors for  a large number of players with mutual interactions, which has led to the theory of mean-field games (MFGs), 
introduced in \cite{laslio07} and \cite{huaetal06}.  

Since the emergence of MFG theory,  the optimal control of mean-field dynamical systems has attracted a lot of interest in the literature, mostly in continuous time. 
It has been first studied in \cite{ahmdin01} by functional analysis method with a value function expressed in terms of the Nisio semigroup of operators.  More recently, several papers 
have adopted the  stochastic maximum principle for characterizing solutions to the controlled McKean-Vlasov systems in terms of adjoint backward stochastic differential equations (BSDEs), 
see \cite{anddje10}, \cite{bucetal11},  \cite{cardel14}. We also refer to  the paper \cite{yon13} which focused on the linear-quadratic (LQ) case where the BSDE from the maximum principle leads to a Riccati equation system. 
It is mentioned in these papers that due to the non-markovian nature of the McKean-Vlasov systems, dynamic programming (also called Bellman optimality) principle does not hold and the problem is time inconsistent in general. 
Indeed, the standard Markov property of the state process, say $X$,  is ruled out, however, as noticed in \cite{benetal13}, this can be restored by working with the marginal distribution of $X$.  The dynamic programming has then been applied independently in  \cite{benetal15} and \cite{laupir14}  for a specific control problem where the objective function depends upon statistics of $X$ like its mean value,  with a mean-field interaction  on the 
drift of the diffusion dynamics of $X$, and in particular  by assuming the existence at all times of a density function for the marginal distribution of $X$.

The purpose of this paper is to provide a detailed analysis of the dynamic programming method  for the optimal control of nonlinear mean-field  systems in discrete time, where the coefficients may depend 
both upon the marginal distributions of the state and of the control.  The case of continuous time McKean-Vlasov equations requires more technicalities and mathematical tools, and will be addressed in \cite{phawei15}. 
The discrete time framework has been also considered in \cite{eletal13} for LQ problem, and arises naturally in situations where signal values are 
available only at certain times.  On the other hand, it can also be viewed as  the discrete time version or approximation of the optimal control of continuous time McKean-Vlasov stochastic differential equations.  
Our methodology is the following.  By u\-sing closed-loop (also called feedback) controls, we first convert the stochastic optimal control problem  into a deterministic control problem involving only the marginal distribution of the state process. 
We then derive the deterministic evolution of the controlled marginal distribution, and prove in its general form the dynamic programming principle (DPP).  This gives sufficient conditions for optimality in terms of calculus of variations 
in the space of feedback control functions. Classical DPP  for stochastic control problem without mean-field interaction falls within our approach. 
We finally apply our method for solving explicitly the mean-variance portfolio selection problem and the multivariate LQ mean-field control problem, and retrieve in particular the results obtained in 
\cite{eletal13} by four diffe\-rent approaches. 

The outline of the paper is as follows.  The next section formulates the McKean-Vlasov control problem in discrete time. In Section \ref{secdynpro}, we develop  the dynamic programming method in this framework. 
Section \ref{secappli} is devoted to applications of the DPP with explicit solutions in the LQ case including the mean-variance problem.

\section{McKean-Vlasov control problem}

\setcounter{equation}{0} \setcounter{Assumption}{0}
\setcounter{Theorem}{0} \setcounter{Proposition}{0}
\setcounter{Corollary}{0} \setcounter{Lemma}{0}
\setcounter{Definition}{0} \setcounter{Remark}{0}

We consider a general class of optimal control of mean-field type in discrete time.  We are given two  measurable spaces $(E,\Bc(E))$ and $(A,\Bc(A))$ representing respectively  
the state space, and  the control space.  We denote by $\Pc(E)$ and $\Pc(A)$  the set of probability measures on $(E,\Bc(E))$ and $(A,\Bc(A))$. 
On  a probability space $(\Omega,\Fc,\P)$, we consider a controlled stochastic dynamics of McKean-Vlasov type: 
\beq \label{mckean}
X_{k+1}^\alpha &=& F_{k+1}(X_k^\alpha,\P_{_{X_k^\alpha}},\alpha_k,\P_{_{\alpha_k}},\eps_{k+1}), \;\;\; k \in \N,  \;\;\; X_0^\alpha \; = \; \xi, 
\enq
for some measurable functions $F_k$ defined from  $E\times\Pc(E)\times A\times\Pc(A)\times\R^d$ into $E$, 
where $(\eps_k)_k$ is a sequence of i.i.d.  random variables, independent of the initial random value $\xi$, and 
we denote by $\F$ $=$ $(\Fc_k)_k$ the filtration with $\Fc_k$ the $\sigma$-algebra generated by $\{\xi,\eps_1,\ldots,\eps_k\}$.  
Here, $(X_k^\alpha)_k$ is the  state process valued in $E$ controlled by the $\F$-adapted process $(\alpha_k)_k$ valued in $A$,  
and we adopted the usual notation in the sequel of the paper: given a random variable $Y$ on $(\Omega,\Fc,\P)$, 
$\P_{_Y}$ denotes its probability distribution under $\P$.  Thus, the dynamics of $(X_k^\alpha)$ depends at any time $k$ of its marginal distribution, 
but also of the marginal distribution of the control, which represents an additional mean-field feature with respect to classical McKean Vlasov equations, and also considered recently  in \cite{eletal13} and \cite{yon13}.

Let us now precise the assumptions on the McKean-Vlasov equation. We shall assume that $(E,|.|)$ is a normed space (most often $\R^d$), 
$(A,|.|)$ is also a normed space (typically a subset of $\R^m$), and we denote by 
$\Pc_{_2}(E)$ the space of square integrable probability measures over $E$ , i.e. $\mu$ $\in$ $\Pc(E)$ s.t. $\|\mu\|_{_2}^2$ $:=$ $\int_E |x|^2 \mu(dx)$ $<$ $\infty$, and 
similarly for $\Pc_{_2}(A)$.  
For any $(x,\mu,a,\lambda)$ $\in$  $E\times\Pc(E)\times A\times\Pc(A)$, and $k$ $\in$ $\N$, we denote by $P_{k+1}(x,\mu,a,\lambda,dx')$ 
the probability distribution of the $E$-valued random variable $F_{k+1}(x,\mu,a,\lambda,\eps_{k+1})$ on $(\Omega,\Fc,\P)$, and we assume

\vspace{2mm}

\noindent 
{\bf (H1)}  For any $k$ $\in$ $\N$, there exists some positive constant $C_{k,F}$ s.t. for all 
$(x,a,\mu,\lambda)$ $\in$ $E\times A\times\Pc(E)\times\Pc(A)$: 
\beqs
\int_E | x'| ^2 P_{k+1}(x,\mu,a,\lambda,dx') & = &  \E\Big[ \big|F_{k+1}(x,\mu,a,\lambda,\eps_{k+1})\big|^2 \Big]  \\ 
& \leq &  C_{k,F} (1 + |x|^2 + |a|^2 + \|\mu\|_{_2}^2 + \|\lambda\|_{_2}^2).
\enqs

Assuming that the initial random value $\xi$ is square integrable, and considering admissible controls $\alpha$ which are square integrable, i.e. 
$\E|\alpha_k|^2$ $<$ $\infty$, for any $k$, it is then clear under {\bf (H1)}  that  $\E |X_k^\alpha|^2$ $<$ $\infty$, i.e. $\P_{_{X_k^\alpha}} \in \Pc_{_2}(E)$, and there exists some positive constant $C_k$ s.t. 
\beq \label{squareX}
\E |X_k^\alpha|^2 & \leq  & C_{k}\big( 1 + \E|\xi|^2 + \sum_{j=0}^{k-1} \E |\alpha_j|^2 \big). 
\enq

The cost functional associated to the system \reff{mckean}  over a finite horizon $n$ $\in$ $\N\setminus\{0\}$ is: 
\beq \label{gainJ}
J(\alpha) & := & \E \Big[ \sum_{k=0}^{n-1} f_k(X_k^\alpha,\P_{_{X_k^\alpha}},\alpha_k,\P_{_{\alpha_k}}) \; + \; g(X_n^\alpha,\P_{_{X_n^\alpha}}) \Big], 
\enq
for any  square integrable $\F$-adapted processes $\alpha$ valued in $A$,  where the running cost functions $f_k$, $k$ $=$ $0,\ldots,n-1$, are measurable real-valued functions on  $E\times\Pc_{_2}(E)\times A\times\Pc_{_2}(A)$, 
and the terminal cost function $g$ is a real-valued measurable function on $E\times\Pc_{_2}(E)$.  
We shall assume 

\vspace{2mm}

\noindent 
{\bf (H2)}  There exist some positive constant $C_{g}$ and for any  $k$ $=$ $0,\ldots,n-1$,
some positive constant $C_{k,f}$  s.t. for all  $(x,a,\mu,\lambda)$ $\in$ $E\times A\times\Pc_{_2}(E)\times\Pc_{_2}(A)$: 
\beqs
 \big|f_k(x,\mu,a,\lambda) \big|  & \leq &  C_{k,f} (1 + |x|^2 + |a|^2 +  \|\mu\|_{_2}^2 + \|\lambda\|_{_2}^2), \\
 \big|g(x,\mu) \big|  & \leq &  C_g (1 + |x|^2 +  \|\mu\|_{_2}^2). 
\enqs

Under {\bf (H1)}-{\bf (H2)}, the cost functional $J(\alpha)$ is well-defined and finite for any admissible control, and  
the objective is to minimize over all admissible controls the cost functional, i.e. by solving
\beq \label{defV0}
V_0 & := & \inf_{\alpha} J(\alpha),
\enq
and when $V_0$ $>$ $-\infty$, find an optimal control $\alpha^*$  i.e. achieving the minimum in \reff{defV0} if it exists.

Problem \reff{mckean}-\reff{defV0} arises in the study of collective behaviors of a large number of players (particles) resulting  from mean-field interactions: typically, the controlled dynamics of a system of $N$ symmetric 
particles are given by 
\beqs
X_{k+1}^{i,\alpha^i} &=& F_{k+1}(X_k^{i,\alpha^i}, \frac{1}{N} \sum_{j=1}^N \delta_{_{X_k^{j,\alpha^i}}},\alpha_k^i,
\frac{1}{N} \sum_{j=1}^N \delta_{_{\alpha_k^{j}}},\eps_{k+1}^i), \;\;\; i=1,\ldots,N,
\enqs
(here $\delta_x$ is the Dirac measure at $x$) and by assuming that a center decides of  the general same policy $\alpha^i$ $=$ $\alpha$ for all  players with same running and terminal gain functions, 
the propagation of chaos argument from McKean-Vlasov theory (see \cite{sni89}) states that when the number of players $N$ goes to infinity, the problem of each agent is asymptotically reduced to the problem  of a single agent 
with controlled dynamics \reff{mckean} and objective \reff{defV0}.  We refer to \cite{cardel13} for a detailed discussion about optimal control of McKean-Vlasov equations and connection with equilibrium of large populations of individuals 
with mean-field interactions.  

\section{Dynamic programming} \label{secdynpro}

\setcounter{equation}{0} \setcounter{Assumption}{0}
\setcounter{Theorem}{0} \setcounter{Proposition}{0}
\setcounter{Corollary}{0} \setcounter{Lemma}{0}
\setcounter{Definition}{0} \setcounter{Remark}{0}

In this section, we make the standing assumptions {\bf (H1)}-{\bf (H2)}, and our purpose is to show  that dynamic programming principle holds 
for problem \reff{defV0},  which we would like to combine with some Markov property of the controlled state process. However, notice that  the McKean-Vlasov type dependence on the dynamics of the state process rules out the standard Markov property  of the controlled process $(X_k^\alpha)$. Actually, this Markov property can be restored by considering 
its  probability law  $(\P_{_{X_k^\alpha}})_k$.  To be more precise and for the sake of definiteness, we shall restrict ourselves to controls 
$\alpha$ $=$ $(\alpha_k)_k$ given in closed-loop (or feedback) form: 
\beq \label{closed}
\alpha_k &=& \tilde\alpha_k(X_k^\alpha), \;\;\;\;\;  k =0,\ldots,n-1, 
\enq
for some deterministic measurable functions $\tilde\alpha_k$ of the state. Notice that the feedback control  may also depend on the (deterministic) marginal distribution, and it will be indeed the case for the optimal one, but to alleviate notation, we omit this dependence which is implicit through  the deterministic function $\tilde\alpha_k$.   
We  denote by $A^E$ the set of measurable functions on $E$ valued in $A$, which satisfy a linear growth condition, and by  $\Ac$ the set of 
admissible controls $\alpha$ in closed loop form \reff{closed} with $\tilde\alpha_k$ in $A^E$, $k$ $\in$ $\N$.  We shall often identify $\alpha$ $\in$ $\Ac$ with the sequence 
$(\tilde\alpha_k)_k$ in $A^E$ via \reff{closed}. Notice that any $\alpha$ $\in$ $\Ac$ satisfies  the square integrability condition, i.e. $\E|\alpha_k|^2$ $<$ $\infty$, for all $k$. Indeed from the linear growth condition on 
$\tilde\alpha_k$ in $A^E$, we have  $\E|\alpha_k|^2$ $\leq$ $C_\alpha(1+ \E|X_k^\alpha|^2)$ for some  constant $C_\alpha$ (depending on $\alpha$),  which gives the square integrability condition by \reff{squareX}.

Next, we show that the initial stochastic control problem can be reduced to a determi\-nistic control problem.  Indeed, the key point is to observe by definition of $\P_{_{X_k^\alpha}}$ and  noting that $\P_{_{\alpha_k}}$ is the image by $\tilde\alpha_k$ of  $\P_{_{X_k^\alpha}}$ for a feedback control $\alpha$ $\in$ $\Ac$,  
that the gain functional in \reff{gainJ} can be rewritten as:   
\beq \label{gainJ2}
J(\alpha) &=& \sum_{k=0}^{n-1} \hat f_k(\P_{_{X_k^\alpha}},\tilde\alpha_k) + \hat g(\P_{_{X_n^\alpha}}), 
\enq
where $\hat f_k$, $k$ $=$ $0,\ldots,n-1$, are defined on $\Pc_{_2}(E)\times A^E$,  $\hat g$ is defined on $\Pc_{_2}(E)$ by:
\beq \label{defhatf}
\hat f_k(\mu,\tilde\alpha) \; := \; \int_E f_k(x,\mu,\tilde\alpha(x),\tilde\alpha\star \mu) \mu(dx),& & \hat g(\mu) \; := \; \int_E g(x,\mu) \mu(dx), 
\enq
and  $\tilde\alpha\star\mu$ $\in$ $\Pc_{_2}(A)$  denotes the image by $\tilde\alpha$ $\in$ $A^E$ of the measure $\mu$ $\in$ $\Pc_{_2}(E)$:
\beqs
\big(\tilde\alpha\star\mu\big)(B)  &=&  \mu\big( \tilde\alpha^{-1}(B) \big), \;\;\; \forall B \in \Bc(A). 
\enqs

Hence, the original problem \reff{defV0} is transformed into a deterministic control pro\-blem involving the infinite dimensional  marginal distribution process. 
Let us then define the dynamic version for  problem \reff{defV0}:  
\beq \label{Vkalpha}
V_k^\alpha & := & \inf_{\beta\in\Ac_k(\alpha)}  \sum_{j=k}^{n-1} \hat f_j(\P_{_{X_j^\beta}},\tilde\beta_j) 
\; + \; \hat g(\P_{_{X_n^\beta}}), \;\;\; k =0,\ldots,n, 
\enq
for $\alpha$ $\in$ $\Ac$, where $\Ac_k(\alpha)$ $=$ $\{\beta\in\Ac: \beta_j = \alpha_j, j=0,\ldots,k-1\}$, with the convention that $\Ac_0(\alpha)$ $=$ $\Ac$, so that $V_0$ $=$ $\inf_{\alpha\in\Ac} J(\alpha)$  is equal to  $V_0^\alpha$.  It is clear that $V_k^\alpha$ $<$ $\infty$, and we shall assume that 
\beq \label{Vkfini}
V_k^\alpha & > & - \infty, \;\;\;\;\;\;\;  k=0,\ldots,n, \;\; \alpha \in \Ac. 
\enq

\begin{Remark}
{\rm The finiteness condition \reff{Vkfini} can be checked a priori directly from the assump\-tions on the model. For example, when $f_k$, $g$, hence $\hat f_k$, $g$, $k$ $=$ $0,\ldots,n-1$, are lower-bounded functions, condition 
\reff{Vkfini} clearly holds.  Another  example is the case when $f_k(x,\mu,a,\lambda)$, $k$ $=$ $0,\ldots,n-1$, and $g$ are lower bounded by a quadratic function  in $x$, $\mu$, and $\lambda$, so that by the 
linear growth condition on $\tilde\alpha$,
\beqs
\hat f_k(\mu,\tilde\alpha) + \hat g(x,\mu) & \geq & - C_k \big( 1 + \|\mu\|_{_2} \big), \;\;\; \forall \mu \in \Pc_{_2}(E), \tilde\alpha\in A^E, 
\enqs
and we are able to derive moment estimates on $X_k^\alpha$, uniformly in $\alpha$: $\big\|\P_{_{X_k^\alpha}}\big\|_{_2}^2$ $=$ $\E[|X_k^\alpha|^2]$ $\leq$ $C_k$, which arises typically when $A$ is bounded from \reff{squareX}. Then, it is clear that \reff{Vkfini} holds true. Otherwise, this finiteness  condition can be checked a posteriori from a verification theorem, see Theorem \ref{Theoverif}. 
}
\ep
\end{Remark}

The  dynamic programming principle (DPP) for the  deterministic control pro\-blem \reff{Vkalpha} takes the following formulation:

\begin{Lemma} (Dynamic Programming Principle)

\noindent Under \reff{Vkfini}, we have
\begin{equation} \label{DPPstandard}
\left\{
\begin{array}{ccl}
V_n^\alpha &=& \hat g(\P_{_{X_n^\alpha}}) \\
V_k^\alpha &=& \Inf_{\beta\in\Ac_k(\alpha)}  \hat f_k(\P_{_{X_k^\beta}},\tilde\beta_k) + V_{k+1}^\beta, \;\;\; k=0,\ldots,n-1. 
\end{array}
\right.
\end{equation}
\end{Lemma}
{\bf Proof.} In the context of deterministic control problem, the proof of the DPP is standard and does not require any measurable selection arguments. For sake of completeness and since it is quite elementary, we give  it.  Denote by 
$J_k(\alpha)$ the cost functional at time $k$, i.e. 
\beqs
J_k(\alpha) &:=& \sum_{j=k}^{n-1} \hat f_k(\P_{_{X_k^\alpha}},\tilde\alpha_k) + \hat g(\P_{_{X_n^\alpha}}), \;\;\; k =0,\ldots,n, 
\enqs
so that $V_k^\alpha$ $=$ $\inf_{\beta\in\Ac_k(\alpha)}J_k(\beta)$ , and by $W_k^\alpha$ the r.h.s. of \reff{DPPstandard}.  Then, 
\beqs
W_k^\alpha & =&  \Inf_{\beta\in\Ac_k(\alpha)} \big[   \hat f_k(\P_{_{X_k^\beta}},\tilde\beta_k) +  \inf_{\gamma\in\Ac_{k+1}(\beta)} J_{k+1}(\gamma) \big] \\
&=& \Inf_{\beta\in\Ac_k(\alpha)} \inf_{\gamma\in\Ac_{k+1}(\beta)} \big[   \hat f_k(\P_{_{X_k^\beta}},\tilde\beta_k) +   J_{k+1}(\gamma) \big] \\
&=& \Inf_{\beta\in\Ac_k(\alpha)} \inf_{\gamma\in\Ac_{k+1}(\beta)} \big[   \hat f_k(\P_{_{X_k^\gamma}},\tilde\gamma_k) +   J_{k+1}(\gamma) \big] \\
&=& \inf_{\gamma\in  \{ \Ac_{k+1}(\beta): \beta \in \Ac_k(\alpha)\}   }  J_{k}(\gamma), 
\enqs
where we used in the third equality the fact that $X_k^\beta$ $=$ $X_k^\gamma$, $\beta_k$ $=$ $\gamma_k$  for $\gamma$ $\in$ $\Ac_{k+1}(\beta)$. 
Finally, we notice that $\{ \Ac_{k+1}(\beta): \beta \in \Ac_k(\alpha)\}$ $=$ $\Ac_k(\alpha)$.  Indeed, the inclusion $\subset$ is clear while for the converse inclusion, it suffices to observe that any 
$\gamma$ in $\Ac_k(\alpha)$ satisfies obviously $\gamma$ $\in$ $\Ac_{k+1}(\gamma)$.  This proves the required equality: $W_k^\alpha$ $=$ $V_k^\alpha$. 
\ep

\vspace{3mm}

Let us now show how one can simplify the DPP by exploiting  the flow  property of  $(\P_{_{X_k^\alpha}})_k$ for any admissible control $\alpha$ in feedback form  $\in$ $\Ac$. Actually, we can derive the evolution of the controlled deterministic process $(\P_{_{X_k^\alpha}})_k$.

\begin{Lemma} \label{fokker}
For any admissible control in closed-loop form $\alpha$ $\in$ $\Ac$, we have
\beq \label{Phimu}
\P_{_{X_{k+1}^\alpha}} &=&  \Phi_{k+1}\big(\P_{_{X_k^\alpha}},\tilde\alpha_k\big),  \;\;\; k \in \N, \;\; \P_{_{X_0^\alpha}} = \P_\xi
\enq
where $\Phi_{k+1}$ is the measurable function defined from $\Pc_{_2}(E)\times A^E$ into $\Pc_{_2}(E)$ by:
\beq \label{defPhi}
\Phi_{k+1}(\mu,\tilde\alpha)(dx') &=& \int_{E} \mu(dx)  P_{k+1}(x,\mu,\tilde\alpha(x),\tilde\alpha\star\mu,dx').  
\enq
\end{Lemma}
{\bf Proof.} Fix $\alpha$ $\in$ $\Ac$. Recall from the definition of the transition probability  
$P_{k+1}(x,\mu,a,\lambda,dx')$ associated to \reff{mckean} that  
\beq \label{Pk}
\P\big[ X_{k+1}^\alpha \in dx'  \big| \Fc_k \big]  &=& P_{k+1}(X_k^\alpha,\P_{_{X_k^\alpha}},\alpha_k,\P_{_{\alpha_k}},dx'), \;\;\; k \in \N.  
\enq
For any bounded measurable function $\varphi$ on $E$, we have by the law of iterated conditional expectation and \reff{Pk}:
\beqs
\E \Big[ \varphi(X_{k+1}^\alpha) \Big] &=& \E \Big[ \E \big[ \varphi(X_{k+1}^\alpha) \big| \Fc_k \big] \Big] \\
&=& \E \Big[ \int_E \varphi(x') P_{k+1} (X_k^\alpha,\P_{_{X_k^\alpha}},\alpha_k,\P_{_{\alpha_k}},dx') \Big] \\
&=& \E \Big[ \int_{E\times E} \varphi(x') P_{k+1} (x,\P_{_{X_k^\alpha}},\tilde\alpha_k(x),\tilde\alpha_k\star \P_{_{X_k^\alpha}},dx') \P_{_{X_k^\alpha}}(dx) \Big] 
\enqs
where we used in the last equality the fact that $\alpha_k$ $=$ $\tilde\alpha_k(X_k^\alpha)$ is in closed loop form,  
the definition of $\P_{_{X_k^\alpha}}$, and noting that $\P_{_{\alpha_k}}$ $=$ $\tilde\alpha_k\star \P_{_{X_k^\alpha}}$. 
This shows the required inductive relation for $\P_{_{X_k^\alpha}}$. 
\ep

\begin{Remark}
{\rm Relation \reff{Phimu} is the Fokker-Planck equation in discrete time for the marginal distribution of the controlled process $(X_k^\alpha)$. In absence of control and McKean-Vlasov type dependence, i.e. $P_{k+1}(x,dx')$ 
does not depend on $(\mu,a,\lambda)$, we retrieve the standard Fokker-Planck equation with a linear function $\Phi_{k+1}(\mu)$ $=$ $\mu P_{k+1}$. In our McKean-Vlasov control context, the function $\Phi_{k+1}(\mu,\tilde\alpha)$  is nonlinear in $\mu$. 
\ep
}
\end{Remark}


 
 \vspace{2mm}

By exploiting  the inductive relation \reff{Phimu} on the  controlled process $(\P_{_{X_k^\alpha}})_k$, the calculation of the value processes $V_k^\alpha$  can be reduced to the  recursive computation of deterministic functions (called value functions) on $\Pc(E)$.

\begin{Theorem} \label{theoDPP} (Dynamic programming and value functions) 

\noindent  Under \reff{Vkfini}, we have  for any $\alpha$ $\in$ $\Ac$, $V_k^\alpha$ $=$ $v_k(\P_{_{X_k^\alpha}})$, $k$ $=$ $0,\ldots,n$,  where $(v_k)_k$ is the sequence of value functions defined recursively on $\Pc_{_2}(E)$ by:
\begin{equation}
\left\{
\begin{array}{ccl} \label{dynpromarkov}
v_n(\mu) &=& \hat g(\mu) \\
v_k(\mu)  &=&  \Inf_{\tilde\alpha \in A^E} \Big[  \hat f_k(\mu,\tilde\alpha) +  v_{k+1}\big(\Phi_{k+1}(\mu,\tilde\alpha)\big)  \Big] 
\end{array}
\right.
\end{equation}
for $k=0,\ldots,n-1$, $\mu$ $\in$ $\Pc_{_2}(E)$.
\end{Theorem}
{\bf Proof.} First observe that for any $\beta$ $\in$ $\Ac_k(\alpha)$, $X_k^\beta$ $=$ $X_k^\alpha$, $k$ $=$ $0,\ldots,n$.  Let us prove the result by 
backward induction. For $k$ $=$ $n$, the result clearly holds since $V_n^\alpha$ $=$ $\hat g(\P_{_{X_n^\alpha}})$. 
Suppose now that at time $k+1$, $V_{k+1}^\alpha$ $=$ $v_{k+1}(\P_{_{X_{k+1}^\alpha}})$ for some deterministic function $v_{k+1}$ and any $\alpha$ $\in$ $\Ac$. Then, from the DPP \reff{DPPstandard} and Lemma \ref{fokker}, we get
\beq 
V_k^\alpha &=& 
 \Inf_{\beta\in\Ac_k(\alpha)}  \hat f_k(\P_{_{X_k^\alpha}},\tilde\beta_k)   + v_{k+1}(\P_{_{X_{k+1}^\beta}})   \nonumber \\
 &=&   \Inf_{\beta\in\Ac_k(\alpha)} w_k(\P_{_{X_k^\alpha}},\tilde\beta_k)  \label{interDPP}
 \enq
 where 
 \beqs
 w_k(\mu,\tilde\alpha) &:= &  \hat f_k(\mu,\tilde\alpha)   +   
   v_{k+1}\big(\Phi_{k+1}\big(\mu,\tilde\alpha\big) \big). 
\enqs
Now, for any $\beta$ $\in$ $\Ac_k(\alpha)$, and since  $\tilde\beta_k$ is valued in $A^E$,  we clearly have: $w_k(\mu,\beta_k)$ 
$\geq$ $\inf_{\tilde\alpha \in A^E} w_k(\mu,\tilde\alpha)$, and so 
$\inf_{\beta\in\Ac_k(\alpha)}w_k(\mu,\tilde\beta_k)$ $\geq$ $\inf_{\tilde\alpha\in A^E} w_k(\mu,\tilde\alpha)$.  Conversely, for any 
$\tilde\alpha$ $\in$ $A^E$, the control $\beta$ defined by $\beta_j$ $=$ $\alpha_j$, $j$ $\leq$ $k-1$, 
and $\tilde\beta_j$ $=$ $\tilde\alpha$ for $j$ $\geq$ $k$, lies in $\Ac_k(\alpha)$,  so: 
$w_k(\mu,\tilde\alpha)$ $\geq$ $\inf_{\beta\in\Ac_k(\alpha)} w_k(\mu,\tilde\beta_k)$, and thus 
$\inf_{\beta\in\Ac_k(\alpha)}w_k(\mu,\tilde\beta_k)$ $=$ $\inf_{\tilde\alpha\in A^E} w_k(\mu,\tilde\alpha)$.  
We conclude from \reff{interDPP} that: $V_k^\alpha$ $=$ 
$v_k(\P_{_{X_k^\alpha}})$ with $v_k(\mu)$ $=$ $\Inf_{\tilde\alpha\in A^E} w_k(\mu,\tilde\alpha)$, i.e.  given by \reff{dynpromarkov}.  
\ep 

\vspace{2mm}

\begin{Remark}
{\rm Problem \reff{defV0}  includes the case where the cost functional in \reff{gainJ} is a nonlinear function of the expected value of the state process, i.e.  the running cost functions and the terminal gain function are in the form:  
$f_k(X_k^\alpha,\P_{_{X_k^\alpha}},\alpha_k)$ $=$ $\bar f_k(X_k^\alpha,\E[X_k^\alpha],\alpha_k)$, $k$ $=$ $0,\ldots,n-1$, $g(X_n^\alpha,\P_{_{X_n^\alpha}})$ $=$ $\bar g(X_n^\alpha,\E[X_n^\alpha])$,  which arise for example in mean-variance pro\-blem (see Section \ref{secappli}).  It is claimed in \cite{bjomur08} and \cite{yon13}  that Bellman optimality principle does not hold, and therefore the problem is time-inconsistent.  This is true when one takes into account only the 
state process $X^\alpha$ (that is its realization), since it is not Markovian, but as shown in this section, dynamic programming principle holds whenever we consider the marginal distribution as state variable.  This gives more information and the price to paid is the infinite-dimensional feature of the marginal distribution state variable.   
}
\ep
\end{Remark}

We complete the above Bellman's optimality principle with a verification theorem, which gives a sufficient condition for finding an optimal control.

\begin{Theorem} \label{Theoverif}
(Verification theorem)

\noindent  (i) Suppose we can find a sequence of real-valued functions $w_k$, $k$ $=$ $0,\ldots,n$, defined on $\Pc_{_2}(E)$ and 
satisfying the dynamic programming relation:
\begin{equation}
\left\{
\begin{array}{ccl} \label{dynpromarkovw}
w_n(\mu) &=& \hat g(\mu) \\
w_k(\mu)  &=&  \Inf_{\tilde\alpha \in A^E} \Big[  \hat f_k(\mu,\tilde\alpha) +  w_{k+1}\big(\Phi_{k+1}(\mu,\tilde\alpha)\big)  \Big] 
\end{array}
\right.
\end{equation}
for $k=0,\ldots,n-1$, $\mu$ $\in$ $\Pc_{_2}(E)$. Then $V_k^\alpha$ $=$ $w_k(\P_{_{X_k^\alpha}})$, 
for all $k$ $=$ $0,\ldots,n$, $\alpha$ $\in$ $\Ac$, and thus $w_k$ $=$ $v_k$.   

\noindent (ii) Moreover,  suppose that at any time $k$ and $\mu$ $\in$ $\Pc(E)$, the infimum in \reff{dynpromarkovw} for $w_k(\mu)$ is attained, by some 
$\tilde\alpha_k^*(.,\mu)$ in $A^E$. Then, by defining by induction the control $\alpha^*$ in feedback form by 
$\alpha_k^*$ $=$ $\tilde\alpha_k^*(X_k^{\alpha^*},\P_{_{X_k^{\alpha^*}}})$, $k$ $=$ $0,\ldots,n-1$, we have  
\beqs
V_0 &=& J(\alpha^*),
\enqs
which means that $\alpha^*$ $\in$ $\Ac$ is an optimal control.
\end{Theorem}
{\bf Proof.} (i) Fix some  $\alpha$ $\in$ $\Ac$, and arbitrary $\beta$ $\in$ $\Ac$ associated to a feedback sequence $(\tilde\beta_k)_k$ in $A^E$. 
Then, from the dynamic programming relation \reff{dynpromarkovw} for $w_k$, and recalling 
the evolution \reff{Phimu} of the controlled marginal distribution $\P_{_{X_k^\beta}}$, we have 
\beqs
w_k(\P_{_{X_k^\beta}}) & \leq & \hat f_k(\P_{_{X_k^\beta}},\tilde\beta_k) + v_{k+1}(\P_{_{X_{k+1}^\beta}}), \;\;\; k = 0,\ldots,n-1.
\enqs
By induction and since $w_n$ $=$ $\hat g$, this gives
\beqs
w_k(\P_{_{X_k^\beta}}) & \leq & \sum_{j=k}^{n-1}  \hat f_j(\P_{_{X_j^\beta}},\tilde\beta_j)  + \hat g(\P_{_{X_n^\beta}}). 
\enqs
By noting that $\P_{_{X_k^\alpha}}$ $=$ $\P_{_{X_k^\beta}}$, when $\beta$ $\in$ $\Ac_k(\alpha)$, and since $\beta$ is arbitrary, this proves that 
$w_k(\P_{_{X_k^\alpha}})$ $\leq$ $V_k^\alpha$.  In particular, $V_k^\alpha$ $>$ $-\infty$, i.e. relation \reff{Vkfini} holds, and then by Theorem 
\ref{theoDPP},  $V_k^\alpha$ is characterized by the sequence of value functions  $(v_k)_k$ defined by the DP \reff{dynpromarkov}. This DP obviously defines by backward induction a unique sequence of functions on $\Pc_{_2}(E)$, hence $w_k$ $=$ $v_k$, $k$ $=$ $0,\dots,n$, and therefore 
$V_k^\alpha$ $=$ $w_k(\P_{_{X_k^\alpha}})$. 

\noindent (ii)  By definition of $\tilde\alpha_k^*$ which attains the infimum in \reff{dynpromarkovw}, we have
\beqs
w_k\big(\P_{_{X_k^{\alpha^*}}}\big) &=&  \hat f_k(\P_{_{X_k^{\alpha^*}}}, \tilde\alpha_k^*(.,\P_{_{X_k^{\alpha^*}}})) + w_{k+1}\big( \P_{_{X_{k+1}^{\alpha^*}}}\big), \;\;\; k =0,\ldots,n-1. 
\enqs
By induction this implies  that 
\beqs
V_0 \; = \; w_0(\P_\xi) &=& \sum_{k=0}^{n-1}  \hat f_k\big(\P_{_{X_k^{\alpha^*}}}, \tilde\alpha_k^*(.,\P_{_{X_k^{\alpha^*}}}) \big) 
\; + \; \hat g(\P_{_{X_n^{\alpha^*}}})  \; = \; J(\alpha^*), 
\enqs
which  shows that $\alpha^*$ is an optimal control. 
\ep

\vspace{2mm}

The above verification theorem, which consists in solving the dynamic programming relation \reff{dynpromarkovw}, is useful to check a posteriori the finiteness condition \reff{Vkfini}, and can be applied in practice fo find explicit solutions to some McKean-Vlasov control problems, as investigated in the next section.

 \section{Applications} \label{secappli}

\setcounter{equation}{0} \setcounter{Assumption}{0}
\setcounter{Theorem}{0} \setcounter{Proposition}{0}
\setcounter{Corollary}{0} \setcounter{Lemma}{0}
\setcounter{Definition}{0} \setcounter{Remark}{0}

\subsection{Special cases}

We consider some particular cases, and provide the special forms of the DPP.

\subsubsection{No mean-field interaction}

We first consider the standard control case where there is no mean-field interaction in the dynamics of the state process,  i.e. 
$F_{k+1}(x,a,\eps_{k+1})$, hence $P_{k+1}(x,a,dx')$ do not depend on $\mu,\lambda$, as well as  in the cost functions 
$f_k(x,a)$ and $g(x)$. For simplicity, we assume that $A$ is a bounded set, which ensures the finiteness condition \reff{Vkfini}. 
In this case, we can see  that  the value functions $v_k$ are in the form
\beq \label{formvkclassi}
v_k(\mu) &=& \int_E \tilde v_k(x) \mu(dx),  \;\;\; k = 0,\ldots,n, 
\enq
where the functions $\tilde v_k$ defined on $E$ satisfy the classical dynamic programming principle:
\begin{equation}
\left\{
\begin{array}{ccl} \label{dynpromarkovclassi}
\tilde v_n(x) &=&  g(x) \\
\tilde v_k(x)  &=&  \Inf_{a \in A} \Big[  f_k(x,a)  +  \E \big[ \tilde v_{k+1} (X_{k+1}^\alpha) \big| X_k^\alpha=x, \alpha_k = a \big]  \Big],
\end{array}
\right.
\end{equation}
for $k$ $=$ $0,\ldots,n-1$.  Let us check this result by backward induction. This holds true for $k$ $=$ $n$ since 
$v_n(\mu)$ $=$ $\hat g(\mu)$ $=$ $\int_E g(x) \mu(dx)$. Suppose that \reff{formvkclassi} holds true at time $k+1$. Then, from the DPP \reff{dynpromarkov}, \reff{defPhi} and  Fubini's theorem, we have
\beqs
v_k(\mu)  &=&  \Inf_{\tilde\alpha \in A^E} \Big[   \int_E f_k(x,\tilde\alpha(x)) \mu(dx)  + 
\int_E  \tilde v_{k+1}(x') \Phi_{k+1}(\mu,\tilde\alpha)(dx')  \Big]  \\
&=& \Inf_{\tilde\alpha \in A^E} \Big[  \int_E \big[ f_k(x,\tilde\alpha(x)) +   \int_E \tilde v_{k+1}(x') P_{k+1}(x,\tilde\alpha(x),dx') \big] \mu(dx) \Big] \\
&= & \Inf_{\tilde\alpha \in A^E}  \int_E \tilde w_k(x,\tilde\alpha(x)) \mu(dx)  
\enqs
where we set $\tilde w_k(x,a)$ $=$ $f_k(x,a) + \int_E \tilde v_{k+1}(x') P_{k+1}(x,a,dx')$. Now, we observe that 
\beq \label{inega}
\Inf_{\tilde\alpha \in A^E}  \int_E \tilde w_k(x,\tilde\alpha(x)) \mu(dx)  &=& \int_E \inf_{a \in A} \tilde w_k(x,a) \mu(dx). 
\enq
Indeed, since for any $\tilde\alpha$ $\in$ $A^E$, the value $\tilde\alpha(x)$ is valued in $A$ for any $x$ $\in$ $E$, it is clear that the inequality $\geq$ in 
\reff{inega} holds true. Conversely, for any $\eps$ $>$ $0$, and $x$ $\in$ $E$, one can find $\tilde\alpha^\eps(x)$ in $A$  such that
\beqs
\tilde w_k(x,\tilde\alpha^\eps(x)) & \leq & \inf_{a\in A} \tilde w_k(x,a) + \eps. 
\enqs
By a measurable selection theorem, the map $x$ $\mapsto$ $\tilde\alpha^\eps(x)$ can be chosen measurable, and since $A$ is bounded, 
the function $\tilde\alpha^\eps$ lies in $A^E$. It follows that 
\beqs
\Inf_{\tilde\alpha \in A^E}  \int_E \tilde w_k(x,\tilde\alpha(x)) \mu(dx) \;\; \leq \;\; \int_E \tilde w_k(x,\tilde\alpha^\eps(x)) \mu(dx) & \leq & \int_E \inf_{a\in A} \tilde w_k(x,a) \mu(dx) + \eps,
\enqs
which shows \reff{inega} since $\eps$ is arbitrary. Therefore, we have $v_k(\mu)$ $=$ $\int_E \tilde v_k(x) \mu(dx)$ with  
\beqs
\tilde v_k(x)  &=&  \inf_{a \in A} \tilde w_k(x,a) \\
&=& \inf_{a \in A} \big[ f_k(x,a) + \int_E \tilde v_{k+1}(x') P_{k+1}(x,a,dx') \big], 
\enqs
which is the relation \reff{dynpromarkovclassi}  at time $k$ from the definition of the transition probability $P_{k+1}$.

\subsubsection{First order interactions}

We consider the case of first order interactions, i.e. the dependence of the model coefficients upon the probability measures is linear in the sense that 
for any $(x,\mu,a)$ $\in$ $E\times\Pc_{_2}(E)\times A$, $\tilde\alpha$ $\in$ $A^E$, 
\beqs
P_{k+1}(x,\mu,a,\tilde\alpha\star\mu,dx') &=& \int_{E} \tilde P_{k+1}(x,y,a,\tilde\alpha(y),dx') \mu(dy), \\
f_k(x,\mu,a,\tilde\alpha\star\mu) & = &  \int_{E}  \tilde f_k(x,y,a,\tilde\alpha(y)) \mu(dy), \;\;\;\;\; g(x,\mu) \; = \; \int_E \tilde g(x,y) \mu(dy),
\enqs
for some  transition probability kernels $\tilde P_{k+1}$ from $E\times E\times A\times A$ into $E$, 
measurable functions   $\tilde f_k$ defined on $E\times E\times A\times A$, $k$ $=$ $0,\ldots,n-1$, and 
$\tilde g$ defined on $E\times E$. 
In this case, the value functions $v_k$ are in the reduced form
\beqs
v_k(\mu) &=& \int_{E^{2^{n-k+1}}} \tilde v_k (\bx_{2^{n-k+1}}) \mu(d\bx_{2^{n-k+1}}), \;\;\; k = 0,\ldots,n,
\enqs
where we denote by  $\bx_{p}$  the $p$-tuple  $(x_1,\ldots,x_p)$ $\in$ $E^p$, by  $\mu(d\bx_p)$  the product measure 
$\mu(dx_1)\otimes\ldots\otimes\mu(dx_p)$, and the functions $\tilde v_k$ are defined recursively on $E^{2^{n-k+1}}$ by 
\begin{equation*}  
\left\{
\begin{array}{ccl} 
\tilde v_n(x,y) &=&  \tilde g(x,y)  \\
\tilde v_k(\bx_{2^{n-k}},\by_{2^{n-k}})  &=&  \Inf_{\tilde\alpha \in A^E} 
\Big[  \tilde f_k(x_1,y_{1},\tilde\alpha(x_{1}),\tilde\alpha(y_1))   \\
& &  \; + \;    \Int_{E^{2^{n-k}}}  \tilde v_{k+1}(\bx_{2^{n-k}}') 
{\bf \tilde P_{k+1}(x_{2^{n-k}},y_{2^{n-k}},\tilde\alpha(x_{2^{n-k}}),\tilde\alpha(y_{2^{n-k}}),dx_{2^{n-k}}') } \Big], 
\end{array}
\right.
\end{equation*}
where we set 
\beqs
& &  {\bf \tilde P_{k+1}(x_{p},y_{p},\tilde\alpha(x_{p}),\tilde\alpha(y_{p}),dx_{p}') } \\
& =  &  \tilde P_{k+1}(x_1,y_1,\tilde\alpha(x_1),\tilde\alpha(y_1),dx_1') \otimes \ldots  \otimes 
\tilde P_{k+1}(x_p,y_p,\tilde\alpha(x_p),\tilde\alpha(y_p),dx_p'). 
\enqs
This result is easily checked  by induction from the DPP \reff{dynpromarkov}, and it is left to  the reader.  

\subsection{Linear-quadratic McKean-Vlasov control problem}

We consider a general multivariate linear McKean-Vlasov dynamics in $E$ $=$ $\R^d$ with control valued in $A$ $=$ $\R^m$: 
\beq 
X_{k+1}^\alpha &=&   \big( B_k X_k^\alpha + \bar B_k \E [X_k^\alpha]  + C_k \alpha_k + \bar C_k \E[\alpha_k] \big)  \label{dynXLQ} \\
& & \;\;\; + \;  \big(  D_k X_k^\alpha + \bar D_k \E[X_k^\alpha] + H_k \alpha_k + \bar H_k \E[\alpha_k])  \eps_{k+1}, \;\;\; k=0,\ldots,n-1,  \nonumber 
\enq
starting from $X_0^\alpha$ $=$ $\xi$, where $B_k$, $\bar B_k$, $D_k$, $\bar D_k$ are constant matrices in $\R^{d\times d}$,   
$C_k$, $\bar C_k$, $H_k$, $\bar H_k$  are constant matrices in $\R^{d\times m}$, and  ($\eps_k)$ is a sequence of i.i.d.  random variables with distribution $\Nc(0,1)$, independent of $\xi$.    
The quadratic cost functional to be minimized is given by
\beq  
J(\alpha) &=& \E \Big[ \sum_{k=0}^{n-1} \big[  (X_k^\alpha)\trans Q_k X_k^\alpha   +   \big(\E[X_k^\alpha]\big)\trans\bar Q_k \E[X_k^\alpha]  
+ L_k\trans X_k^{\alpha}   + \bar L_k\trans \E[X_k^{\alpha}]   \nonumber \\
& &\hspace{1.5cm} + \;   \alpha_k\trans R_k\alpha_k + \big(\E[\alpha_k]\big)\trans \bar R_k\E[\alpha_k] \big] \;  \nonumber \\
& & \hspace{1.5cm}  + \;    \big(X_n^\alpha\big)\trans Q X_n^\alpha +  \big(\E[X_n^\alpha]\big)\trans\bar Q  \E[X_n^\alpha] 
+ L\trans X_n^{\alpha}   + \bar L\trans \E[X_k^{\alpha}]   \Big], \label{JLQ}
\enq
 for some constants matrices $Q_k$, $\bar Q_k$, $Q$, $\bar Q$,  in $\R^{d\times d}$, $R_k$, $\bar R_k$ in $\R^{m\times m}$,  
 and vectors $L_k$, $\bar L_k$, $L$, $\bar L$ $\in$ $\R^d$,  $k$ $=$ $0,\ldots,n-1$.  Here 
 $x\trans$ denotes the transpose of a matrix/vector $x$. 
 Since the cost functions are real-valued,  we may assume w.l.o.g. that all these matrices $Q_k$, $\bar Q_k$, $Q$, $\bar Q$, $R_k$ and $\bar R_k$  are symmetric. 
This model is in the form \reff{mckean} and  associated to a transition probability  satisfying: 
\beq 
P_{k+1}(x,\mu,a,\lambda,dx') & \leadsto & \Nc\Big(M_k(x,\mu,a,\lambda) ; \Sigma_k(x,\mu,a,\lambda)\Sigma_k(x,\mu,a,\lambda)\trans\Big) \label{loiPk} \\
M_k(x,\mu,a,\lambda) &=& B_k x + \bar B_k \bar\mu + C_k a + \bar C_k \bar\lambda \nonumber \\
\Sigma_k(x,\mu,a,\lambda) &=& D_k x + \bar D_k \bar\mu + H_k a + \bar H_k \bar\lambda \nonumber 
\enq
where we set for any $\mu$ $\in$ $\Pc_2(\R^d)$ (resp. $\Pc_{_2}(\R^m)$)  symmetric matrix $\Lambda$ $\in$ $\R^{d\times d}$ (resp. in $\R^{m\times m}$):
\beqs
\bar\mu \; := \; \int_{} x \mu(dx),    \;\;\; 
\bar\mu_{_2}(\Lambda)  \; :=  \;  \int_{} x\trans \Lambda x \mu(dx),  \;\;\;    {\rm Var}(\mu)(\Lambda) \; := \; \bar\mu_{_2}(\Lambda)  - \bar\mu\trans\Lambda\bar\mu,  
\enqs
and in  the form  \reff{gainJ}, hence \reff{gainJ2} for feedback controls, with 
\beqs
\hat f_k(\mu,\tilde\alpha) &=&  {\rm Var}(\mu)(Q_k)    +     \bar\mu\trans(Q_k+\bar Q_k)\bar\mu  +  (L_k + \bar L_k)\trans\bar\mu  \\
& &  \;\;\; + \; {\rm Var}(\tilde\alpha\star\mu)(R_k) +  \overline{\tilde\alpha\star\mu}\trans(R_k+\bar R_k)\overline{\tilde\alpha\star\mu} \\
\hat g(\mu) &=&  {\rm Var}(\mu)(Q)  \; + \;     \bar\mu\trans(Q+\bar Q)\bar\mu \; + \; (L + \bar L)\trans\bar\mu.  
\enqs

We look for  candidate  $w_k$, $k$ $=$ $0,\ldots,n$, of values functions  satisfying the dynamic programming principle \reff{dynpromarkov},  in the quadratic form: 
\beq \label{formvk2}
w_k(\mu) &=&  {\rm Var}(\mu)(\Lambda_k)  +  \bar\mu\trans\Gamma_k\bar\mu + \rho_k\trans\bar\mu + \chi_k,
\enq
for some constant symmetric matrices $\Lambda_k$ and  $\Gamma_k$ in $\R^{d\times d}$, vector $\rho_k$ $\in$ $\R^d$ and real 
$\chi_k$  to be determined below.   
We proceed  by backward induction. For $k$ $=$ $n$, we see that $w_k$ $=$ $\hat g$ ($=$ $v_k$) iff
\beq \label{termLQ}
\Lambda_n \; = \; Q, & &  \Gamma_n \; = \;  Q + \bar Q, \;\; \rho_n \; = \; L + \bar L, \; \chi_n \; = \; 0. 
\enq
Now, suppose that the form \reff{formvk2} holds true at time $k+1$, and observe from \reff{defPhi} and  \reff{loiPk} that for any $\mu$ $\in$ $\Pc_2(\R^d)$, $\tilde\alpha$ $\in$ 
$A^E$, $\Lambda$ $\in$ $\R^{d\times d}$, we have by Fubini's theorem: 
\beqs
\overline{\Phi_{k+1}(\mu,\tilde\alpha)} &=&  \int_{\R^d} \E\big[ Y(x,\mu,\tilde\alpha) \big]  \mu(dx) \\
\overline{\Phi_{k+1}(\mu,\tilde\alpha)}_{_{2}}(\Lambda) &=&  \int_{\R^d} \E\big[ Y(x,\mu,\tilde\alpha)\trans\Lambda Y(x,\mu,\tilde\alpha) \big]  \mu(dx),
\enqs
where $Y(x,\mu,\tilde\alpha)$ $\leadsto$  $\Nc\Big(M_k(x,\mu,\tilde\alpha(x),\tilde\alpha\star\mu) ; 
\Sigma_k(x,\mu,\tilde\alpha(x),\tilde\alpha\star\mu)\Sigma_k(x,\mu,\tilde\alpha(x),\tilde\alpha\star\mu)\trans\Big)$. Therefore, 
\beqs
\overline{\Phi_{k+1}(\mu,\tilde\alpha)}  &=& (B_k + \bar B_k) \bar\mu  + (C_k + \bar C_k)  \overline{\tilde\alpha\star\mu},
\enqs 
and after some tedious but straightforward calculation:
\beqs
{\rm Var}(\Phi_{k+1}(\mu,\tilde\alpha))(\Lambda) &=& \overline{\Phi_{k+1}(\mu,\tilde\alpha)}_{_{2}}(\Lambda) -  \overline{\Phi_{k+1}(\mu,\tilde\alpha)}\trans\Lambda\overline{\Phi_{k+1}(\mu,\tilde\alpha)} \\
&=&   \int_{\R^d} 
\big[ \Sigma_k(x,\mu,\tilde\alpha(x),\tilde\alpha\star\mu)\trans \Lambda \Sigma(x,\mu,\tilde\alpha(x),\tilde\alpha\star\mu)   \\
& & \;\;\;\;\;\;\; + \;\; M_k(x,\mu,\tilde\alpha(x),\tilde\alpha\star\mu)\trans \Lambda M(x,\mu,\tilde\alpha(x),\tilde\alpha\star\mu) \big] \mu(dx)  \\
& & \; - \;   \Big( (B_k + \bar B_k) \bar\mu  + (C_k + \bar C_k)  \overline{\tilde\alpha\star\mu} \Big)\trans\Lambda \Big( (B_k + \bar B_k) \bar\mu  + (C_k + \bar C_k)  \overline{\tilde\alpha\star\mu} \Big) \\
&=& {\rm Var}(\mu)( B_k\trans\Lambda B_k +D_k\trans\Lambda D_k) +  \bar\mu\trans(D_k +\bar D_k)\trans \Lambda (D_k + \bar D_k) \bar\mu \\
& & \; \; + \; {\rm Var}(\tilde\alpha * \mu)(H_k\trans \Lambda H_k  + C_k\trans \Lambda C_k) \\
& & \;\; + \;   \overline{\tilde\alpha\star\mu}\trans (H_k +\bar H_k) \trans \Lambda (H_k + \bar H_k)  \overline{\tilde\alpha\star\mu}  \\
& & \;\; + \;\;  2\int_{\R^d} (x-\bar\mu) \trans (D_k\trans\Lambda H_k + B_k\trans \Lambda C_k)\tilde\alpha(x) \mu(dx) \\
& & \;\;+\;\;   2\bar\mu\trans  (D_k + \bar D_k) \trans \Lambda (H_k + \bar H_k)   \overline{\tilde\alpha\star\mu}.  
\enqs 
Then, $w_k$ satisfies the DPP \reff{dynpromarkov} iff
\beq 
w_k(\mu) &=&    \inf_{\tilde\alpha\in A^E} \Big[  \hat f_k(\mu,\tilde\alpha) +  {\rm Var}(\Phi_{k+1}(\mu,\tilde\alpha))(\Lambda_{k+1})  + 
\overline{\Phi_{k+1}(\mu,\tilde\alpha)}\trans\Gamma_{k+1}\overline{\Phi_{k+1}(\mu,\tilde\alpha)} \Big]  \label{wkLQ} \\
&=&   {\rm Var}(\mu)(Q_k +  B_k\trans\Lambda_{k+1} B_k + D_k\trans\Lambda_{k+1} D_k) +  \inf_{\tilde\alpha\in A^E}  G_{k+1}^\mu(\tilde\alpha)  \nonumber \\
& & \; + \; \bar\mu\trans\big(Q_k + \bar Q_k +   (D_k + \bar D_k)\trans\Lambda_{k+1}(D_k + \bar D_k) +   (B_k + \bar B_k)\trans\Gamma_{k+1}(B_k + \bar B_k) \big) \bar\mu   \nonumber \\
& & \; + \; \big(L_k + \bar L_k + (B_k + \bar B_k)\trans\rho_{k+1} \big)\trans\bar\mu  + \chi_{k+1},
\enq
 where we define  the function $G_{k+1}^\mu$ $:$ $L^2(\mu;A)$ $\mapsto$ $\R$ by
\beq
G_{k+1}^\mu(\tilde\alpha) &=& {\rm Var}(\tilde\alpha\star\mu)(V_k)  \;+\;  \overline{\tilde\alpha\star\mu}\trans W_k   \overline{\tilde\alpha\star\mu} 
 \;+\; 2\int_{\R^d} (x-\bar\mu) \trans S_k \tilde\alpha(x) \mu(dx) \nonumber \\
 & & \;+\; 2\bar\mu\trans  T_k \overline{\tilde\alpha\star\mu}  \; + \;  \rho_{k+1}\trans(C_k +\bar C_k)\overline{\tilde\alpha\star\mu},   \label{defGk} 
\enq
and we set $V_k$ $=$ $V_k(\Lambda_{k+1})$, $W_k$ $=$ $W_k(\Lambda_{k+1},\Gamma_{k+1})$, $S_k$ $=$ $S_k(\Lambda_{k+1})$, $T_k$ $=$ $T_k(\Lambda_{k+1},\Gamma_{k+1})$, with 
\begin{equation} \label{defVWST}
\left\{
\begin{array}{ccl}
V_k(\Lambda_{k+1}) &=&R_k + H_k\trans \Lambda_{k+1} H_k + C_k\trans \Lambda_{k+1} C_k;\\
W_k(\Lambda_{k+1},\Gamma_{k+1}) &=& R_k + \bar R_k + (C_k+\bar C_k) \trans \Gamma_{k+1} (C_k + \bar C_k) +(H_k + \bar H_k) \trans \Lambda_{k+1} (H_k + \bar H_k) \\
S_k(\Lambda_{k+1}) &=&D_k\trans\Lambda_{k+1} H_k + B_k\trans\Lambda_{k+1}  C_k ;\\
T_k(\Lambda_{k+1},\Gamma_{k+1}) &=& (D_k+\bar D_k) \trans \Lambda_{k+1} (H_k + \bar H_k ) + (B_k + \bar B_k))\trans \Gamma_{k+1} (C_k +\bar C_k).
\end{array}
\right.
\end{equation} 
Here,  $L^2(\mu;A)$ $\supset$ $A^E$ is   the Hilbert space of measurable functions on $E$ $=$ $\R^d$ valued in $A$ $=$ $\R^m$ and square integrable w.r.t. $\mu$ $\in$ $\Pc_{_2}(E)$.

We now search for the infimum of the function $G_{k+1}^\mu$, and shall make the following assumptions on the  symmetric matrices of the quadratic cost functional and on the coefficients of the state dynamics: 

\vspace{1mm}

\hspace{-7mm} {\bf (c0)} 
\begin{equation} \nonumber 
\left\{
\begin{array}{ccl}
 Q \; \geq \; 0, \; Q+\bar Q \; \geq \; 0, & & Q_k \; \geq \; 0, \; Q_k + \bar Q_k \; \geq \; 0,  \\
 & & R_k \; \geq   \; 0, \; R_k + \bar R_k \; \geq  \; 0, \;\; k=0,\ldots,n-1, 
\end{array}
\right.
\end{equation}
and for all $k$ $=$ $0,\ldots,n-1$ (with the convention that $Q_n$ $=$ $Q$, $\bar Q_n$ $=$ $\bar Q$)
\begin{itemize}
\item[{\bf (c1)}] $R_k > 0$, $\;$ or $\;$  [rank($C_k$) $=$ $d$, $Q_{k+1}$ $>$ $0$], $\;$ or $\;$  [rank($H_k$) $=$ $d$, $Q_{k+1}$ $>$ $0$], 
\item[{\bf (c2)}]  $R_k + \bar R_k> 0$,  or $\;$  [rank($C_k+\bar C_k$) $=$ $d$, $Q_{k+1}+\bar Q_{k+1}$ $>$ $0$], $\;$ or $\;$  [rank($H_k+\bar H_k$) $=$ $d$, $Q_{k+1}$ $>$ $0$].
\end{itemize}
Conditions {\bf (c0)}-{\bf (c1)}-{\bf (c2)} is slightly weaker than the condition in  \cite{eletal13} (see their Theorem 3.1), where the condition {\bf (c0)} is strengthened  to $R_k$ $>$ $0$ and  $R_k + \bar R_k$ $>$ $0$ for all $k$ $=$ $0,\ldots,n-1$, for ensuring the existence of an optimal control.  
We relax this positivity condition  with the conditions {\bf (c1)}-{\bf (c2)} in order to include the case of mean-variance problem (see the example at the end of this section).   Actually, as we shall see in Remark \ref{remconvex}, these conditions  
will guarantee that for   $\Lambda_{k}$, $\Gamma_k$ to be determined below,  the function $G_{k+1}^\mu$ is convex and coercive on  $L^2(\mu;A)$ for any $k$ $=$ $0,\ldots,n-1$.  
For the moment, we derive after some  straightforward calculation  the Gateaux derivative of $G_{k+1}^\mu$ at $\tilde\alpha$ in the direction $\beta$ $\in$ $L^2(\mu;A)$, which is given by:
\beqs
DG_{k+1}^\mu(\tilde\alpha ; \beta)  &:=& \lim_{\eps\rightarrow 0}  \frac{G_{k+1}^\mu(\tilde\alpha + \eps\beta) - G_{k+1}^\mu(\tilde\alpha)}{\eps}
\; = \; \int_{\R^d} g_{k+1}(x,\tilde\alpha) \beta(x) \mu(dx)
\enqs
with
\beqs
g_{k+1}(x, \tilde\alpha)&=& 2\tilde\alpha(x)\trans V_k  + 2 \overline{\tilde\alpha\star\mu}  \trans  \big(W_k -V_k\big) \\
& & \;\;\; + \;  2(x-\mu)\trans S_k  +  2 \bar\mu\trans T_k + \rho_{k+1}\trans(C_k + \bar C_k).
\enqs
We shall check later in Remark \ref{remconvex} that $V_k$ and $W_k$ are positive symmetric matrices, hence invertible.  We thus  see that $DG_{k+1}^\mu(\tilde\alpha;.)$ vanishes for $\tilde\alpha$ $=$  $\tilde\alpha_k^*(.,\mu)$ s.t. 
$g_{k+1}(x,\tilde\alpha_k^*(.,\mu))$ $=$ $0$ for all $x$ $\in$ $\R^d$, which gives:
\beq
\tilde\alpha_k^*(x,\mu)  &=&   - V_k^{-1} S_k\trans(x-\bar\mu)\;  -\;   W_k^{-1}T_k \trans \bar\mu
\; - \;  \frac{1}{2} W_k^{-1}(C_k + \bar C_k) \trans \rho_{k+1}
 \label{optialphaLQ} 
\enq
and then  after some straightforward caculation:
\beqs
G_{k+1}^\mu(\tilde\alpha_k^*(.,\mu))&= & - {\rm Var}(\mu)\big(  S_k  V_k^{-1} S_k \trans \big) \;  - \;  
\bar\mu\trans \big( T_k W_k^{-1}T_k\trans\big)\bar\mu  \;- \; \bar\mu \trans T_k W_k^{-1} (C_k + \bar C_k) \trans \rho_{k+1}  \\
& & \; -\;  \frac{1}{4}  \rho_{k+1}\trans(C_k + \bar C_k)W_k^{-1}(C_k + \bar C_k) \trans \rho_{k+1}. 
\enqs
Assuming for the moment that $\tilde\alpha_k^*(.,\mu)$ attains the infimum of $G_{k+1}^\mu$ (this is a consequence of the convexity and coercivity of  $G_{k+1}^\mu$ shown in Remark \ref{remconvex}), and plugging the above expression in \reff{wkLQ},  we see that 
$w_k$ is like the function $\mu$ $\mapsto$ $G_{k+1}^\mu(\tilde\alpha_k^*(.,\mu))$,  a linear combination of terms in ${\rm Var(\mu)}(.)$, 
$\bar\mu\trans(.)\bar\mu$,  and by identification with the form \reff{formvk2}, we obtain an inductive relation for $\Lambda_k$, $\Gamma_k$, $\rho_k$, 
$\chi_k$:  
\begin{equation} \label{Riccati}
\left\{
\begin{array}{ccl}
\Lambda_k &=&  Q_k + B_k \trans \Lambda_{k+1} B_k +  D_k \trans \Lambda_{k+1} D_k -S_k(\Lambda_{k+1}) V_k^{-1}(\Lambda_{k+1})  
S_k\trans(\Lambda_{k+1}) \\ 
\Gamma_k  &=& (Q_k +\bar Q_k) + (B_k + \bar B_k) \trans \Gamma_{k+1} (B_k + \bar B_k) + (D_k + \bar D_k) \trans \Lambda_{k+1} (D_k + \bar D_k)\\
& & \; - \; T_k(\Lambda_{k+1},\Gamma_{k+1}) W_k^{-1}(\Lambda_{k+1},\Gamma_{k+1}) T_k\trans(\Lambda_{k+1},\Gamma_{k+1}) \\
\rho_k &=& L_k + \bar L_k +\big[(B_k + \bar B_k) - (C_k +\bar C_k) W_k^{-1}(\Lambda_{k+1},\Gamma_{k+1}) 
T_k\trans(\Lambda_{k+1},\Gamma_{k+1}) \big]\rho_{k+1} \\
\chi_k &=& \chi_{k+1} - \frac{1}{4}  \rho_{k+1}\trans(C_k + \bar C_k) W_k^{-1}(\Lambda_{k+1},\Gamma_{k+1})(C_k + \bar C_k) \trans \rho_{k+1}.
\end{array}
\right.
\end{equation}
for all $k$ $=$ $0,\ldots,n-1$, starting from the terminal condition \reff{termLQ}.  The relations for $(\Lambda_k,\Gamma_k)$ in \reff{Riccati} are 
two algebraic  Riccati difference equations, while the equations for $\rho_k$ and $\chi_k$  are linear equations once 
$(\Lambda_k,\Gamma_k)$ are determined.  This system \reff{Riccati} is the same as the one obtained in \cite{eletal13}. 
In the particular mean-variance problem considered at the end of this section, we can obtain explicit closed-form expressions for the solutions 
$(\Lambda_k,\Gamma_k,\rho_k,\chi_k)$ to this Riccati system. However, in general, there are no closed-form formulae, and these quantities  are simply computed by induction.

 \vspace{1mm} 
  
In the  following remark, we check the issues that have left open up to now.

\begin{Remark} \label{remconvex}
{\rm  Let  conditions {\bf (c0)}-{\bf (c1)}-{\bf (c2)} hold. We prove by backward induction that  for all $k$ $=$ $1,\ldots,n$, the matrices  $V_{k-1}$ $=$ $V_{k-1}(\Lambda_{k})$, $W_{k-1}$ $=$ $W_{k-1}(\Lambda_{k},\Gamma_{k})$ are 
symmetric positive, hence invertible, with $(\Lambda_k,\Gamma_k)$ given by \reff{Riccati}, together with the nonnegativity of the symmetric matrices 
$\Lambda_k$, $\Gamma_k$, which will immediately gives the convexity and coercivity of 
the function $G_{k}^\mu$ in \reff{defGk} for $\mu$ $\in$ $\Pc_{_2}(\R^d)$. 

At time $k$ $=$ $n$, we have $\Lambda_n$ $=$ $Q$ $\geq$ $0$,  $\Gamma_n$ $=$ $Q+\bar Q$ $\geq$ $0$, and thus from \reff{defVWST},  the $d\times d$ matrices 
$V_{n-1}$ $=$ $V_{n-1}(\Lambda_{n})$, $W_{n-1}$ $=$ $W_{n-1}(\Lambda_{n},\Gamma_{n})$ are symmetric  positive under {\bf (c0)}-{\bf (c1)}-{\bf (c2)}: indeed, for example if $R_{n-1}$ $=$ $0$, then the rank condition 
on $C_{n-1}$ or $H_{n-1}$ together with the positivity of $\Lambda_n$ $=$ $Q$ in {\bf (c1)} will ensure that $V_{n-1}$ is positive. Now, suppose that the assertion is true at time $k+1$, i.e. 
$V_k$, $W_k$ are symmetric positive, and $\Lambda_{k+1}$, $\Gamma_{k+1}$ are symmetric nonnegative.  Then, it is clear from \reff{Riccati} that $\Lambda_k$ and $\Gamma_k$ are symmetric,  and 
noting that they can be rewritten from the expression of $V_k,W_k,S_k,T_k$ in \reff{defVWST} as
\begin{equation} \nonumber 
\left\{
\begin{array}{ccl}
\Lambda_k &=& Q_k +  S_k V_k^{-1} R_k  \big(S_k V_k^{-1}\big) \trans  \; + \;  \Big(B_k - C_k  (S_k V_k^{-1})\trans \big]\trans \Lambda_{k+1} \big[B_k - C_k (S_k V_k^{-1})\trans \Big) \\ 
& & \; + \;  \Big( D_k - H_k (S_kV_k^{-1})\trans \Big)\trans \Lambda_{k+1} \Big(D_k - H_k (S_k V_k^{-1}) \trans \Big) \\
\Gamma_k &=& Q_k +\bar Q_k \; + \;   T_k W_k^{-1}(R_k +\bar R_k) \big(T_k W_k^{-1}\big)\trans \\
& & \; + \;   \Big( B_k + \bar B_k  -  (C_k +\bar C_k)(T_k W_k^{-1})\trans \Big)\trans\Gamma_{k+1} \Big( B_k + \bar B_k  -  (C_k +\bar C_k)(T_k W_k^{-1})\trans\Big) \\
& & \; + \;  \Big( D_k +\bar D_k  - (H_k +\bar H_k)(T_k W_k^{-1})\trans \Big)\trans
\Gamma_{k+1}\Big( D_k +\bar D_k  -   (H_k +\bar H_k)(T_k W_k^{-1})\trans \Big),
\end{array}
\right.
\end{equation}
it is also clear that  they  are nonnegative under   {\bf (c0)}. Finally from the expression \reff{defVWST} at time $k-1$, we see that  $V_{k-1}$ $=$ $V_{k-1}(\Lambda_k)$ and $W_{k-1}$ $=$ $W_{k-1}(\Lambda_k,\Gamma_k)$ 
are symmetric positive under {\bf (c0)}-{\bf (c1)}-{\bf (c2)}, which shows the required assertion. 
}
\ep
\end{Remark}

\vspace{1mm}

In view of the above derivation and Remark \ref{remconvex}, it follows that the functions $w_k$, $k$ $=$ $0,\ldots,n$, given in the quadratic form \reff{formvk2} with $(\Lambda_k,\Gamma_k,\rho_k,\chi_k)$ as in \reff{Riccati}, satisfy  by construction the 
DPP \reff{dynpromarkov}, and by the verification theorem, this implies that the value functions are given by $v_k$ $=$ $w_k$, while the optimal control is given in feedback form from \reff{optialphaLQ} by:
\beq \label{alphafeedback}
\alpha_k^* &=&  \tilde\alpha_k(X_k^{*},\P_{_{X_k^{*}}}) \; = \; 
- V_k^{-1} S_k\trans\big(X_k^* - \E [ X_k^*] \big)\;  -\;   W_k^{-1}T_k \trans \E [ X_k^* ],
\enq
where  $X_k^*$ $=$ $X_k^{\alpha^*}$ is the optimal wealth process with the feedback control $\alpha^*$.  We retrieve the expression obtained in \cite{eletal13}  (see e.g. their Theorem 3.1). 
We can push further our calculations  to get an explicit form of the optimal control expressed only in terms of  the state process 
(and not on its mean).  Indeed, from the linear dynamics \reff{dynXLQ}, we have
\beqs
\E[X_{k+1}^*] & =& (B_k + \bar B_k) \mathbb{E} [X_k^*] + (C_k + \bar C_k) \mathbb{E} [\alpha_k^*]\\
& =& (B_k + \bar B_k) \mathbb{E} [X_k^*]  - (C_k + \bar C_k)   \big(W_k\trans\big)^{-1}T_k \trans \E [ X_k^* ] \; = \;  N_k \mathbb{E} [X_k^*],
\enqs
with  $N_k= B_k + \bar B_k -(C_k + \bar C_k) W_k^{-1}T_k$, for $k$ $=$ $0,\ldots,n-1$, and so by induction:
\beqs
\E [X_k^*] &=& N_{k-1}\ldots N_0  \E[\xi].  
\enqs 
Plugging into \reff{alphafeedback}, this gives the explicit form of the optimal control as
\beq \label{optialphanonfeedback}
\alpha_k^* &=&  - V_k^{-1} S_k\trans X_k^*  \; + \; \big(V_k^{-1}S_k\trans  - W_k^{-1}T_k\trans\big) N_{k-1}\ldots N_0 \E[\xi], \;\;\; k=0,\ldots,n-1. 
\enq
We observe that the optimal control at any time $k$ does not only depend on the current state  $X_k^*$ but also on its the initial  state $\xi$ (via its mean).


\vspace{3mm}

\noindent {\bf Example: Mean-variance portfolio selection}
 
 \vspace{1mm}
 
\noindent The  mean-variance discrete-time problem consists in minimizing the cost functional: 
 \beqs
 J(\alpha) &=& \frac{\gamma}{2} {\rm Var}(X_n^\alpha) -  \E[X_n^\alpha]  \\
 &=& \E\big[ \frac{\gamma}{2} \big(X_n^\alpha\big)^2 -  X_n^\alpha\big] \; - \;  \frac{\gamma}{2} \Big( \E[X_n^\alpha] \Big)^2,  
 \enqs
 for some $\gamma$ $>$ $0$, with a dynamics for the wealth process $(X_k^\alpha)$ valued in $E$ $=$ $\R$  controlled by the amount $\alpha_k$ valued in $A$ $=$ $\R$ 
 invested in the stock at time $k$ (we assume zero interest rate): 
 \beq \label{wealth}
 X_{k+1}^\alpha &=& X_k^\alpha + \alpha_k ( b \Delta + \sigma \sqrt{\Delta} \eps_{k+1}), \;\;\; k =0,\ldots,n-1, \; X_0^\alpha = x_0. 
 \enq
 Here $x_0$ $\in$ $\R$ is the initial capital, $b$, $\sigma$ $>$ $0$ are some constants, representing respectively the rate of return and volatility of the stock, $\Delta$ $>$ $0$ is a parameter, e.g. $\Delta$ $=$ $T/n$, arising when considering a time discretization of  a continuous-time model over $[0,T]$, and ($\eps_k)$ is a sequence of i.i.d.  random variables with distribution $\Nc(0,1)$.  This univariate model fits into the LQ framework 
 \reff{dynXLQ}-\reff{JLQ} with $d$ $=$ $m$ $=$ $1$, and  
 \beqs
 B_k = 1, \; \bar B_k = 0, \; C_k = b\Delta, \; \bar C_k = 0, \; D_k = \bar D_k = 0, \; H_k = \sigma\sqrt{\Delta}, \; \bar H_k = 0,  \\
 Q_k = \bar Q_k = L_k = \bar L_k = R_k = \bar R_k = 0, \;\; Q = \frac{\gamma}{2},  \; \bar Q = - \frac{\gamma}{2}, \; L = 0, \; \bar L = -1.  
 \enqs 
Conditions {\bf (c0)}-{\bf (c1)}-{\bf (c2)} are clearly satisfied, and the Riccati system \reff{Riccati} for $(\Lambda_k,\Gamma_k,\rho_k,\chi_k)$ 
$\in$ $\R_+\times\R_+\times\R\times\R$ is written in this case as: 
\begin{equation} \nonumber
\left\{
\begin{array}{ccl}
\Lambda_k &=& \Lambda_{k+1} \frac{\sigma^2}{\sigma^2 + b^2 \Delta}    \\
\Gamma_k &=&  \frac{\sigma^2 \Lambda_{k+1}}{\sigma^2 \Lambda_{k+1} + b^2\Delta \Gamma_{k+1}}  \Gamma_{k+1}   \\
\rho_k &=& \frac{\sigma^2 \Lambda_{k+1}}{b^2 \Delta \Gamma_{k+1} +  \sigma^2 \Lambda_{k+1}} \rho_{k+1}    \\
\chi_k &=& \chi_{k+1} -    \; \frac{1}{4} \frac{b^2\Delta \rho_{k+1}^2}{\sigma^2 \Lambda_{k+1} + b^2\Delta \Gamma_{k+1}},
\end{array}
\right.
\end{equation}
together with the terminal condition $\Lambda_n$ $=$ $\frac{\gamma}{2}$, $\Gamma_n$ $=$ $0$,  $\rho_n$ $=$ $-1$, $\chi_n$ $=$ $0$. 
This leads by induction to the explicit form for $(\Lambda_k,\Gamma_k,\rho_k,\chi_k)$:
\begin{equation} \label{eqphi}
\left\{
\begin{array}{ccl}
\Lambda_k &=&  \frac{\gamma}{2}  \Big(\frac{\sigma^2}{\sigma^2 + b^2 \Delta} \Big)^{n-k}, \\
\Gamma_k & =& 0, \;\; \rho_k \; = \; -1 \\
\chi_k &=& - \frac{1}{2\gamma} \Big( \Big(\frac{\sigma^2 + b^2 \Delta} {\sigma^2}\Big)^{n-k} - 1 \Big). 
\end{array}
\right.
\end{equation}
The value functions are  then explicitly given by
\beqs
v_k(\mu) &=&  \frac{\gamma}{2}  \Big(\frac{\sigma^2}{\sigma^2 + b^2 \Delta} \Big)^{n-k} {\rm Var(\mu)}  - \bar\mu  - \frac{1}{2\gamma} \Big( \Big(\frac{\sigma^2 + b^2 \Delta} {\sigma^2}\Big)^{n-k} - 1 \Big),
\enqs
for all $k$ $=$ $0,\ldots,n$, $\mu$ $\in$ $\Pc_{_2}(\R)$. Moreover, the optimal control is given in feedback form  from \reff{alphafeedback} by:
\beqs \label{optialpha2}
\alpha_k^* &=& \tilde\alpha_k(X_k^{*},\P_{_{X_k^{*}}}) \; = \;  - \frac{  b }{\sigma^2 + b^2 \Delta} \big( X_k^* -\E[X_k^*]\big) \; + \; \frac{b}{\sigma^2\gamma}  \Big(\frac{\sigma^2 + b^2 \Delta} {\sigma^2}\Big)^{n-k-1}, 
\enqs
where $X_k^*$ $=$ $X_k^{\alpha^*}$ is the optimal wealth process with the feedback control $\alpha^*$. It is then explicitly written from \reff{optialphanonfeedback} by   
\beq \label{optialphafin}
\alpha_k^* &=& - \frac{  b }{\sigma^2 + b^2 \Delta} \Big[ X_k^* - x_0  -  \frac{1}{\gamma}  \Big(1 + \frac{b^2}{\sigma^2}\Delta\Big)^{n} \Big]. 
\enq
 We then observe that the optimal control at any time $k$ does not only depend on the current wealth $X_k^*$ but also on the initial wealth $x_0$.  This expression \reff{optialphafin} of the optimal control is the discrete time 
 analog of the continuous time optimal control  obtained in \cite{lizho00} or \cite{anddje10}. Actually, if we view \reff{wealth} as a time discretization (with a time step $\Delta$ $=$ $T/n$) 
 of a continuous time Black-Scholes model for the stock price over $[0,T]$, with  a controlled wealth dynamics
\beqs
dX_t^\alpha &=& \alpha_t (bdt + \sigma dW_t), \;\;\; X_0^\alpha = x_0,
\enqs
then by sending $n$ to infinity (hence $\Delta$ to zero) into \reff{optialphafin}, we retrieve  the closed-form expression of the optimal control in \cite{lizho00} or \cite{anddje10}: 
 \beqs
 \alpha_t^* &=& - \frac{  b }{\sigma^2} \Big[X_t^{\alpha^*} - x_0 -  \frac{1}{\gamma} \exp\big(\frac{b^2}{\sigma^2}T\big) \Big]. 
 \enqs

\end{document}